\documentclass[12pt]{amsart}
\usepackage{amssymb}
\usepackage[margin=0.755in]{geometry}

\usepackage[hyphenbreaks]{breakurl}
\usepackage[hyphens]{url}
\newcommand{\nc}{\newcommand}
\nc{\issuenumber}{40}
\nc{\issuemonth}{Dec}
\nc{\issueyear}{'16}
\nc{\set}[2]{\{\, #1 : #2\,\}}
\nc{\st}{\op{star}}
\nc{\tPa}[6]{\bibitem{#1} {#2}, \emph{#3}, {#4}, to appear (#5 pages).}
\nc{\sPa}[5]{\bibitem{#1} {#2}, \emph{#3}, {#4}, submitted.}
\nc{\FinSeqs}[1]{{{#1}^{<\alephes}}}
\nc{\rest}[1]{\mid_{#1}}
\nc{\Bgp}{{\Z^\N}}
\nc{\Arh}{Arhangel'ski\u{\i}}
\nc{\grbl}{{\mbox{\textit{\tiny gp}}}}
\long\def\forget#1\forgotten{}

\nc{\Ga}{\Gamma}
\nc{\sub}{\subseteq}
\nc{\R}{\mathbb{R}}
\nc{\Cp}{\mathrm{C}_\mathrm{p}}
\nc{\Op}{\mathrm{O}}
\nc{\Fin}{{[\N]^{<\aleph_0}}}
\nc{\alephes}{{\aleph_0}}
\nc{\ed}{

\general\end{document}}

\nc{\fbx}[1]{\fbox{$#1$}}\nc{\nop}{$\times$}\nc{\fbn}{\!\!\fbox{\!\nop\!}\!\!}
\nc{\yup}{\checkmark}\nc{\fby}{\!\!\fbox{\!\yup\!}\!\!}\nc{\mbq}{\mb{?}}
\nc{\mb}[1]{{\mbox{\textbf{#1}}}}\nc{\smb}[1]{{\!\!\mb{#1}\!\!}}\nc{\x}{\times}
\nc{\<}{\left <}\nc{\Cantor}{{\{0,1\}^\N}}\nc{\oo}{\infty}
\nc{\NR}{{\bbR^\N}}\nc{\Iff}{\Leftrightarrow}\nc{\mypar}[1]{\par\medskip\noindent\textbf{#1.}}
\nc{\roth}{{[\N]^{\alephes}}}\nc{\sr}[2]{{\txt{$#1$\\$#2$}}}
\nc{\gpbl}{{\mbox{\textit{\tiny gp}}}}\nc{\fx}{\mathfrak{x}}\nc{\fb}{\mathfrak{b}}
\nc{\fg}{\mathfrak{g}}\nc{\fc}{\mathfrak{c}}\nc{\fd}{\mathfrak{d}}
\nc{\fp}{\mathfrak{p}}\nc{\fs}{\mathfrak{s}}\nc{\ADD}{{\mathsf   {ADD}}}
\nc{\COV}{{\mathsf   {COV}}}\nc{\NON}{{\mathsf   {NON}}}\nc{\COF}{{\mathsf   {COF}}}
\nc{\sseq}[1]{\{#1 : n\in\N\}}\nc{\Impl}{\Rightarrow}
\nc{\upannouncement}[1]{[\S\ref{#1} above]}\nc{\dnannouncement}[1]{[\S\ref{#1} below]}
\nc{\E}{\exists}\nc{\cI}{\mathcal{I}}\nc{\cN}{\mathcal{N}}\nc{\cP}{\mathcal{P}}
\nc{\cA}{\mathcal{A}}\nc{\cM}{\mathcal{M}}\nc{\Null}{\mathcal{N}}
\nc{\op}{\operatorname}\nc{\cov}{\mathsf{cov}}\nc{\add}{\mathsf{add}}
\nc{\cof}{\mathsf{cof}}\nc{\cf}{\mathsf{cf}}\nc{\non}{\mathsf{non}}\nc{\spst}{\supseteq}
\nc{\CH}{the Continuum Hypothesis}\nc{\bbR}{\mathbb{R}}\nc{\Q}{\mathbb{Q}}
\nc{\EdNote}[1]{\par\medskip\noindent\textbf{#1.}}\nc{\fo}{\mathfrak{od}}
\nc{\cl}[1]{\overline{#1}}\nc{\impl}{\rightarrow}\nc{\arrays}{{{\{0,1\}}^{\N\x\N}}}
\nc{\w}{\omega}\nc{\ft}{\mathfrak{t}}\nc{\h}{\mathfrak{h}}\nc{\Cite}[2]{{\cite[#1]{#2}}}
\nc{\bq}{\begin{quote}}\nc{\eq}{\end{quote}}
\nc{\cK}{\mathcal{K}}\nc{\cB}{\mathcal{B}}\nc{\BG}{\cB_\Gamma}
\nc{\BL}{\cB_\Lambda}\nc{\BT}{\cB_\Tau}\nc{\BTstar}{\cB_{\Tau^*}}\nc{\BO}{\cB_\Omega}
\nc{\CG}{C_\Gamma}\nc{\CL}{C_\Lambda}\nc{\CT}{C_\Tau}\nc{\CTstar}{C_{\Tau^*}}
\nc{\CO}{C_\Omega}\nc{\sone}{\mathsf{S}_1}\nc{\sfin}{\mathsf{S}_\mathrm{fin}}
\nc{\Sc}{\mathsf{S}_c}\nc{\ufin}{\mathsf{U}_\mathrm{fin}}\nc{\gone}{\mathsf{G}_1} \nc{\gfin}{\mathsf{G}_\mathrm{fin}}\nc{\seq}[1]{\{#1\}_{n\in\N}}\nc{\Un}{\bigcup}
\nc{\nin}{\not\in}\nc{\cF}{\mathcal{F}}\nc{\cG}{\mathcal{G}}\nc{\cU}{\mathcal{U}}
\nc{\cV}{\mathcal{V}}\nc{\cW}{\mathcal{W}}\nc{\fU}{\mathfrak{U}}\nc{\fu}{\mathfrak{u}}
\nc{\fV}{\mathfrak{V}}\nc{\fW}{\mathfrak{W}}\nc{\psin}{pseudo-intersection}
\nc{\NN}{{\N^\N}}\nc{\N}{\mathbb{N}}\nc{\bbN}{\mathbb{N}}\nc{\Z}{\mathbb{Z}}
\nc{\as}{\subseteq^*}\nc{\sm}{\setminus}\nc{\sbst}{\subseteq}
\nc{\by}[2]{\par\hfill\emph{#1}, #2}\nc{\nby}[1]{\par\hfill\emph{#1}}\nc{\Tau}{\mathrm{T}}
\nc{\CE}{\textsc{CE}}
\newtheorem{thm}{Theorem}[section]\nc{\bthm}{\begin{thm}} \nc{\ethm}{\end{thm}}
\newtheorem{prop}[thm]{Proposition}\nc{\bprp}{\begin{prop}} \nc{\eprp}{\end{prop}}
\newtheorem{fact}[thm]{Fact}\nc{\bfct}{\begin{fact}} \nc{\efct}{\end{fact}}
\newtheorem{prob}[thm]{Problem}\nc{\bprb}{\begin{prob}} \nc{\eprb}{\end{prob}}
\newtheorem{lem}[thm]{Lemma}\nc{\blem}{\begin{lem}} \nc{\elem}{\end{lem}}
\newtheorem{claim}[thm]{Claim}\nc{\bclm}{\begin{claim}} \nc{\eclm}{\end{claim}}
\newtheorem{cor}[thm]{Corollary}\nc{\bcor}{\begin{cor}} \nc{\ecor}{\end{cor}}
\newtheorem{conj}[thm]{Conjecture}\nc{\bcnj}{\begin{conj}} \nc{\ecnj}{\end{conj}}
\theoremstyle{definition}\newtheorem{defn}[thm]{Definition}\nc{\bdfn}{\begin{defn}} \nc{\edfn}{\end{defn}}
\theoremstyle{remark}\newtheorem{rem}[thm]{Remark}\nc{\brem}{\begin{rem}} \nc{\erem}{\end{rem}}
\newtheorem{cnv}[thm]{Convention}\nc{\bcnv}{\begin{cnv}} \nc{\ecnv}{\end{cnv}}
\newtheorem{exam}[thm]{Example}\nc{\bexm}{\begin{exam}} \nc{\eexm}{\end{exam}}
\nc{\bpf}{\begin{proof}} \nc{\epf}{\end{proof}}
\nc{\be}{\begin{enumerate}}\nc{\ee}{\end{enumerate}}\nc{\bi}{\begin{itemize}}
\nc{\ei}{\end{itemize}}\nc{\itm}{\item}
\nc{\general}{\small\vfill\par\noindent\hrulefill\par
\noindent\textbf{Previous issues.} 
\url{http://front.math.ucdavis.edu/search?\&t=\%22SPM+Bulletin\%22}
\\[0.1cm]
\textbf{Contributions and free subscription.} Email \url{tsaban@math.biu.ac.il}.
}

\nc{\link}[1]{\par\hfill{\url{#1}}}
\nc{\fillin}{{\Huge To be completed}}
\nc{\arXivl}[4]{\subsection{#2}{#4}\par\hfill{\arx{#1}}\par\hfill\emph{#3}}
\nc{\arXiv}[3]{\subsection{#2} {#3}, \arx{#1}\par\hfill}
\nc{\DOIpaper}[5]{\subsection{#2}{#4}\par\hfill{\url{http://dx.doi.org/#1}}\par\hfill\emph{#3}}
\nc{\AMSPaper}[5]{\subsection{#3}{#5}\par\hfill{\url{#1}}\par\hfill\emph{#4}\par\hfill{#2}}
\nc{\nAMSPaper}[4]{\subsection{#2} {#3}, {#4}, \url{#1}}
\nc{\AMS}[3]{\subsection{#1} {#2}, \url{#3}}
\nc{\SPMBul}{\textbf{$\mathcal{SPM}$ Bulletin}}
\nc{\BulEnd}{\par\bigskip\noindent
Boaz Tsaban\\
\emph{E-mail}: tsaban@math.biu.ac.il\\
\emph{URL}: http://www.cs.biu.ac.il/\~{}tsaban}

\nc{\arx}[1]{\url{http://arxiv.org/abs/#1}}

\nc{\probissue}{\emph{Problem of the issue}}

\title[$\mathcal{SPM}$ Bulletin \textbf{\issuenumber} (\issuemonth{} \issueyear)]{%
$\mathcal{SPM}$ Bulletin\\[0.5cm]
Issue \issuenumber{} (\issuemonth{} \issueyear{})}

\begin{document}
\maketitle


\section{Editor's note}

In addition to the beautiful research announcements, this issue includes an announcement of a festive conference dedicated to Selection Principles, celebrating the
60th birthday of Marion Scheepers. This will be an opportunity for all active researchers in the field to exchange new methods and results in the field, and for researchers looking for new, promising research directions to consider this topic in a thorough manner.

\medskip

With best regards,

\by{Boaz Tsaban}{tsaban@math.biu.ac.il}

\hfill \texttt{http://www.cs.biu.ac.il/\~{}tsaban}

\section{New conference: Frontiers of Selection Principles}

Dear Friends and Colleagues,

\medskip

We are glad to invite you to the conference Frontiers of Selection Principles, celebrating the 60th birthday of Marion Scheepers. This festive conference will take place at the Dewajtis Campus, Bielanski forest, Warsaw, 21.8--1.9, 2017: 21.8--25.8 tutorial, 26.8 excursion, 27.8--1.9 conference. (This is immediately after the Logic Colloquium and immediately before the Bedlewo meeting on set theoretic topology and analysis.)

\emph{Selection principles} connect topology, set theory, and functional analysis and make it possible to transport and apply methods from each of these fields to the other ones. It is now one of the most active streams of research within set theory and general topology. This conference will be fully dedicated to selection principles and their applications. It will begin with a one week tutorial, aimed mainly at students and researchers with no prior knowledge of selection principles. This tutorial will provide an overview of the filed, and a detailed introduction to its main methods. Participants with no prior knowledge attending the tutorial will thus fully benefit from the lectures of the second week, and will be able to consider possible directions for research in the field.

The second conference week will consist of invited lectures, by leading experts, on results in the frontiers of selection principles. We recommend that students and participants with expertise in other disciplines attend both weeks of the conference.

The expected accommodation cost is approximately 20EUR per night (full board: 30EUR).

Registration fee: Early registration: 50EUR. Normal registration: 80EUR. Registration fees will help reducing the costs for some students or early career researchers that would otherwise not be able to attend this meeting. Discount possibilities will be provided in future emails.

\medskip

\noindent Tentative list of conference speakers:

\begin{center}
\begin{tabular}{lll}
Leandro Aurichi	&	Liljana Babinkostova	&	Taras Banakh\\
Angello Bella*	&	Daniel Bernal Santos	&	Maddalena Bonanzinga*\\
Lev Bukovsky	&	Steven Clontz*		&	Samuel Da-Silva*\\
Rodrigo Dias		&	Rafa\l{} Filipow			& 	David Gauld\\
Ondrej Kalenda	&	Ljubi\v{s}a Ko\v{c}inac*		&	Adam Krawczyk\\
Adam Kwela		& 	Andrzej Nowik		&	Selma\"O{}zca\v{g}*\\
Yinhe Peng		& 	Szymon Plewik		&	Robert Rałowski	\\		
Masami Sakai		&	Marion Scheepers		&  	Paul Szeptycki*	\\
Piotr Szewczak	&	Franklin Tall			& 	Secil Tokg\"oz	\\
Boaz Tsaban		& 	Tomasz Weiss		& 	Lyubomyr Zdomskyy\\
Shuguo Zhnag	&	Ondrej Zindulka		&	Szymon Zeberski
\end{tabular}
\end{center}

\noindent (* to be confirmed)

\medskip

Subject to time constraints, students and other participants may have an opportunity to contribute a short lecture on a topic of selection principles.

\noindent\textbf{Special request:} Since a part of this workshop is especially accessible to students, we would appreciate your forwarding this message to your students and to other students who may be interested in attending this conference.

\noindent The conference is organized by Cardinal Stefan Wyszyński University in Warsaw.

\noindent\textbf{Important:} This is the only message distributed widely. Please email Piotr Szewczak (p.szewczak@wp.pl) to be included in the mailing list for details and updates.

\medskip

On behalf of the organizing committee,

\medskip

\noindent Piotr Szewczak (Cardinal Stefan Wyszynski University and Bar-Ilan University)\\
\noindent Boaz Tsaban (Bar-Ilan University)\\
\noindent Lyubomyr Zdomskyy (Kurt Gödel Research Center)

\section{Long announcements}

\arXivl{1607.04688}
{On the definability of Menger spaces which are not $\sigma$-compact}
{Franklin D. Tall and Secil Tokgoz}
{Hurewicz proved completely metrizable Menger spaces are $\sigma$-compact. We
extend this to \v{C}ech-complete Menger spaces and consistently to projective
Menger metrizable spaces. On the other hand, it is consistent that there is a
co-analytic Menger space that is not $\sigma$-compact.}

\arXivl{1607.04781}
{Definable versions of Menger's conjecture}
{Franklin D. Tall}
{Menger's conjecture that Menger spaces are $\sigma$-compact is false; it is
true for analytic subspaces of Polish spaces and undecidable for more complex
definable subspaces of Polish spaces. For non-metrizable spaces, analytic
Menger spaces are $\sigma$-compact, but Menger continuous images of co-analytic
spaces need not be. The general co-analytic case is still open, but many
special cases are undecidable, in particular, Menger topological groups. We
also prove that if there is a Michael space, then productively Lindelof
Cech-complete spaces are $\sigma$-compact. We also give numerous
characterizations of proper K-Lusin spaces. Our methods include the Axiom of
Co-analytic Determinacy, non-metrizable descriptive set theory, and
Arhangel'skii's work on generalized metric spaces.}

\arXivl{1608.03546}
{Discrete subsets in topological groups and countable extremally
  disconnected groups}
{Evgenii Reznichenko and Ol'ga Sipacheva}
{It is proved that any countable topological group in which the filter of
neighborhoods of the identity element is not rapid contains a discrete set with
precisely one onisolated point and that the existence of a nondiscrete
countable extremally disconnected group implies the existence of a rapid
ultrafilter.}

\arXivl{1608.06210}
{$\sigma$-Ideals and outer measures on the real line}
{S. Garc\'ia-Ferreira, A. H. Tomita, and Y. F. Ortiz-Castillo}
{A {\it weak selection} on $\mathbb{R}$ is a function $f: [\mathbb{R}]^2 \to
\mathbb{R}$ such that $f(\{x,y\}) \in \{x,y\}$ for each $\{x,y\} \in
[\mathbb{R}]^2$. In this article, we continue with the study (which was
initiated in [ag]) of the outer measures $\lambda_f$ on the real line
$\mathbb{R}$ defined by weak selections $f$. One of the main results is to show
that $CH$ is equivalent to the existence of a weak selection $f$ for which: \[
\mathcal \lambda_f(A)= \begin{cases} 0 & \text{if $|A| \leq \omega$,}\\ \infty
& \text{otherwise.} \end{cases} \] Some conditions are given for a
$\sigma$-ideal of $\mathbb{R}$ in order to be exactly the family
$\mathcal{N}_f$ of $\lambda_f$-null subsets for some weak selection $f$. It is
shown that there are $2^\mathfrak{c}$ pairwise distinct ideals on $\mathbb{R}$
of the form $\mathcal{N}_f$, where $f$ is a weak selection. Also we prove that
Martin Axiom implies the existence of a weak selection $f$ such that
$\mathcal{N}_f$ is exactly the $\sigma$-ideal of meager subsets of
$\mathbb{R}$. Finally, we shall study pairs of weak selections which are
"almost equal" but they have different families of $\lambda_f$-measurable sets.}

\arXivl{1609.05822}
{Selectively (a)-spaces from almost disjoint families are necessarily countable under a certain parametrized weak diamond principle}
{Charles J.G. Morgan and Samuel G. Da Silva}
{The second author has recently shown [20] that any selectively (a) almost disjoint 
family must have cardinality strictly less than $2^{\aleph_0}$, so under the Continuum 
Hypothesis such a family is necessarily countable. However, it is also shown in the 
same paper that $2^{\aleph_0} < 2^{\aleph_1}$ alone does not avoid the existence of 
uncountable selectively (a) almost disjoint families. We show in this paper that a 
certain effective parametrized weak diamond principle is enough to ensure countability 
of the almost disjoint family in this context. We also discuss the deductive strength 
of this specific weak diamond principle (which is consistent with the negation of the 
Continuum Hypothesis, apart from other features).
\par
Houston Journal of Mathematics 42 (2016), 1031--1046.}

\arXivl{1206.0722}
{Notes on the od-Lindel\"of property}
{Mathieu Baillif}
{A space is od-compact (resp. od-Lindel\"of) provided any cover by open dense
sets has a finite (resp. countable) subcover. We first show with simple
examples that these properties behave quite poorly under finite or countable
unions. We then investigate the relations between Lindel\"ofness,
od-Lindel\"ofness and linear Lindel\"ofness (and similar relations with
`compact'). We prove in particular that if a $T_1$ space is od-compact, then
the subset of its non-isolated points is compact. If a $T_1$ space is
od-Lindel\"of, we only get that the subset of its non-isolated points is
linearly Lindel\"of. Though, Lindel\"ofness follows if the space is moreover
locally openly Lindel\"of (i.e. each point has an open Lindel\"of
neighborhood).}

\arXiv{1610.04800}
{Relating games of Menger, countable fan tightness, and selective
  separability}
{Steven Clontz}
{By adapting techniques of Arhangel'skii, Barman, and Dow, we may equate the
existence of perfect-information, Markov, and tactical strategies between two
interesting selection games. These results shed some light on Gruenhage's
question asking whether all strategically selectively separable spaces are
Markov selectively separable.}

\arXivl{1611.04998}
{Convergence in topological groups and the Cohen reals}
{Alexander Shibakov}
{We show that it is consistent to have an uncountable sequential group of
intermediate sequential order while no countable such groups exist. This is
proved by adding $\omega_2$ Cohen reals to a model of $\diamondsuit$.}


\section{Short announcements}\label{RA}

\arXiv{1607.04756}
{Hereditarily normal manifolds of dimension $> 1$ may all be metrizable}
{Alan Dow and Franklin D. Tall}

\arXiv{1607.07188}
{A long chain of P-points}
{Borisa Kuzeljevic and Dilip Raghavan}

\arXiv{1607.07669}
{The cofinal structure of precompact and compact sets in general metric
  spaces}
{Aviv Eshed, M. Vincenta Ferrer, Salvador Hern\'andez, Piotr Szewczak,
  Boaz Tsaban}

\arXiv{1607.07978}
{$\omega^\omega$-bases in topological and uniform spaces}
{Taras Banakh}

\arXiv{1608.03381}
{$I$-convergence classes of sequences and nets in topological spaces}
{Amar Kumar Banerjee and Apurba Banerjee}

\arXiv{1608.07144}
{The existence of continuous weak selections and orderability-type
  properties in products and filter spaces}
{Koichi Motooka, Dmitri Shakhmatov, Takamitsu Yamauchi}

\arXiv{1610.04506}
{Weakly linearly Lindel\"of monotonically normal spaces are Lindel\"of}
{I. Juh\'asz, V. V. Tkachuk, R. G. Wilson}

\arXiv{1611.07267}
{On the cardinality of almost discretely Lindel\"of spaces}
{Santi Spadaro}

\arXiv{1611.08289}
{Categorical properties on the hyperspace of nontrivial convergent
sequences}
{S. Garcia-Ferreira, R. Rojas-Hernandez, Y. F. Ortiz-Castillo}

\arXiv{1612.06651}
{First countable and almost discretely Lindel\"of $T_3$ spaces have
  cardinality at most continuum}
{Istv\'an Juh\'asz, Lajos Soukup, Zolt\'an Szentmikl\'ossy}
  
\ed